\documentclass{sapm}
\usepackage{lscape}  
\usepackage{CJK}     
\usepackage{amssymb} 

\def\today{Submitted 22~August 2015, revised 10~December~2015, accepted 14 January 2016} 

\usepackage{graphicx}
\jname{STUDIES IN APPLIED MATHEMATICS}
\jvol{AA}
\jyear{2016}
\doi{10.1146/((please add article doi))}

\def \ccomma{\raise 2pt\hbox{,}} 
\def \D {\hbox{d}}

\def \Rcubec {\mathbb{R}^3(c)}

\def \Pn     {{\rm Pn}}
\def \PVI    {{\rm P6}}
\def \PV     {{\rm P5}}
\def \PIII   {{\rm P3}}
\def \PI     {{\rm P1}}

\def \coth{\mathop{\rm coth}\nolimits}

\def \directsum {\oplus}

\def \diag {\mathop{\rm diag}\nolimits}
\def \tr   {\mathop{\rm tr}\nolimits}
\def \mod#1{\vert #1 \vert}

\def \barQ {\overline{Q}}

\def \bfF {{\bf F}}
\def \bfN {{\bf N}}

\def \bfsigma {{\bf \sigma}}
\def \barz {{\bar z}}

\def \ch {c_{\rm h}}
\def \cu {c_{\rm u}}
\def \cz {c_{\rm z}}
\def \cq {c_{\rm q}}

\def \Hu{\textit{\rm H}}  
\def \Hqp{\textit{\sl H}} 

\def \as{a_{\rm s}}
\def \ad{a_{\rm d}}

\def \tildev{\tilde v}
\def \tildeh{\tilde h}
\def \tildeq{\tilde q}  
\def \tilder{\tilde r}

\begin{document}

\title{Reductions of Gauss-Codazzi equations}

\author[Robert Conte]
{Robert Conte ${}^{1,2}$
and A.~Michel Grundland ${}^{3,4}$\thanks{ 
    Centre de math\'ematiques et de leurs applications,  
 \'Ecole normale sup\'erieure de Cachan, 
 61, avenue du Pr\'esident Wilson, F--94235 Cachan Cedex, France.
Robert.Conte@cea.fr, Grundlan@crm.umontreal.ca
} \\ \today}
\affil{1. 
    Centre de math\'ematiques et de leurs applications  
\\  \'Ecole normale sup\'erieure de Cachan, CNRS, Universit\'e Paris-Saclay,
\\ 61, avenue du Pr\'esident Wilson, F--94235 Cachan Cedex, France.
}
\affil{2. Department of Mathematics, The University of Hong Kong,
\\ Pokfulam Road, Hong Kong.}

\affil{3. Centre de recherches math\'ematiques, Universit\'e de Montr\'eal
\\ Case postale 6128, Succursale Centre ville,
\\ Montr\'eal, Qu\'ebec H3C 3J7, Canada}

\affil{4. D\'epartement de math\'ematiques et d'informatique,
\\ Universit\'e du Qu\'ebec \`a Trois-Rivi\`eres
\\ Case postale 500, Trois-Rivi\`eres, Qu\'ebec G9A 5H7, Canada}

\maketitle

\begin{abstract}
We prove that conformally parametrized
surfaces in Euclidean space $\Rcubec$ of curvature $c$
admit a symmetry reduction of their Gauss-Codazzi equations 
whose general solution is expressed with the sixth Painlev\'e function. 
Moreover,
it is shown that the two known solutions of this type (Bonnet 1867, 
Bobenko, Eitner and Kitaev 1997) can be recovered by such a reduction.
\end{abstract}


												
02.30.Hq Ordinary differential equations

02.30.Ik Integrable systems 

02.30.Jr Partial differential equations


02.30.+g   
						
02.40.-k Differential geometry

02.40.Hw Classical differential geometry 

20.40.Dr


Primary 53A05 

Secondary 58F07

\vfill\eject
\tableofcontents
\vfill\eject

\section{Introduction}

Consider a three-dimensional Riemannian manifold $\mathbb{R}^3(\kappa)$ 
having a constant curvature $\kappa$.
When $\kappa$ is respectively negative, zero, positive,
this three-dimensional manifold is respectively
the hyperbolic space $\mathbb{H}^3(\kappa)$, the Euclidean space $\mathbb{R}^3$, 
the sphere $\mathbb{S}^3(\kappa)$ of 
radius $\kappa^{-1/2}$.
For convenience we will denote $\kappa=-c^2 \in \mathbb{R}$ 
and also denote the just mentioned spaces $\Rcubec$ as, respectively, 
$\mathbb{H}^3(c), \mathbb{R}^3, \mathbb{S}^3(c)$. 

Let $\bfF$ be an immersion of some complex two-dimensional Riemannian manifold 
into $\Rcubec$.
In conformal coordinates $z=x+i y$, $\barz=x-i y$,
the two fundamental forms of the surface are given by 
\begin{eqnarray}
& &
{\rm I}=<\D \bfF,\D \bfF>=e^u \D z \ \D \barz, 
\label{eqform1} 
\end{eqnarray} 
and
\begin{eqnarray}
& &
{\rm II}=-<\D\bfF,\D\bfN>=Q \D z^2 + e^u H \D z \ \D\barz + \barQ \D\barz^2,
\label{eqform2} 
\end{eqnarray} 
in which 
$\bfN$ is a unit vector normal to the tangent plane,
$u$ and $H$           are real    valued functions,
and 
$Q$ is a complex valued function.

In terms of the principal curvatures $1/R_1$ and $1/R_2$, 
the mean curvature $H$ 
and the total (or Gaussian) curvature $K$ are defined as
\begin{eqnarray}
& & H=\frac{1}{2}\left(\frac{1}{R_1}+\frac{1}{R_2}\right), 
\\
& & K=\frac{1}{R_1 R_2}=H^2 - c^2 - 4 e^{-2 u} \mod{Q}^2
     = -2 e^{-u} u_{z \barz}.
\label{eqdefcurvatures}
\end{eqnarray}

The moving frame defined by
\begin{eqnarray}
& &
\bfsigma=
\left\lbrace 
\begin{array}{ll}
\displaystyle{
{}^{\rm t} (\bfF_z,\bfF_\barz,\bfN) {\hskip 4.0truemm} (c=0),\ 
}\\ \displaystyle{
{}^{\rm t} (\bfF,\bfF_z,\bfF_\barz,\bfN)\ (c\not=0), 
}
\end{array}
\right.
\label{eqMoving-frame} 
\end{eqnarray} 
evolves linearly as 
\begin{eqnarray}
& & {\hskip -8.0truemm}
\bfsigma_z=U \bfsigma,\ \bfsigma_\barz=V \bfsigma,\ 
\label{eqgradsigma}
\end{eqnarray}
in which the two linear operators only depend on $u,Q,\barQ,H,c$.

Instead of representing $U$ and $V$ as third ($c=0$) or fourth ($c\not=0$) order matrices
(or even fifth order ones \cite{Bianchi1903} when $c=0$),
it proves quite convenient to use 
the representation by second order matrices \cite{B1994,Springborn}, 
\begin{eqnarray}
& & {\hskip -22.0truemm}
U=\pmatrix{
 (1/4) u_z           & -Q e^{-u/2} \cr 
 (1/2) (H+c) e^{u/2} & -(1/4) u_z \cr },\
V=\pmatrix{
-(1/4) u_\barz & -(1/2) (H-c) e^{u/2} \cr 
\barQ e^{-u/2} & (1/4) u_\barz \cr}.
\label{eqGW-su2-R3R4}
\end{eqnarray}
The zero-curvature (Maurer-Cartan) condition
\begin{eqnarray}
& & {\hskip -8.0truemm}
[\partial_z-U,\partial_\barz-V]=U_\barz - V_z + U V - V U=0,
\label{eqzerocurvature}
\end{eqnarray}
generates the Gauss-Codazzi system of equations,
\begin{eqnarray}
& &
\left\lbrace 
\begin{array}{ll}
\displaystyle{
u_{z \barz} + \frac{1}{2} (H^2-c^2) e^u -2 \mod{Q}^2 e^{-u}=0\ \hbox{ (Gauss)},
}\\ \displaystyle{
Q_\barz-\frac{1}{2} H_z e^u =0,\ 
\barQ_z-\frac{1}{2} H_\barz e^u =0\ \hbox{(Codazzi)}.
}
\end{array}
\right.
\label{eqGaussCodazziR3c}
\end{eqnarray}

After the initial work of the Polish school \cite{Sym1982,Sym1985,CGS1995},
who interpreted, among others, 
isothermic surfaces in $\mathbb{R}^3$ as ``soliton surfaces'',
the relation between geometry and integrability is now well established, 
see the lecture notes \cite{B1999} for a summary.
Many integrable 1+1-dimensional partial differential equations (PDEs)
are contained 
 in the Gauss-Codazzi system (\ref{eqGaussCodazziR3c}) \cite{CGS1995,ChenLi2002}.
In particular, two solutions have been found,
separated by more than one century \cite{Bonnet1867,BEK1997},
in terms of the master Painlev\'e function \PVI\ (Appendix \ref{AppendixPVI}) 
or its confluence to \PV.
Other solutions have also been found in terms of \PIII\ \cite[Eq.~(20)]{BEK1997}.
\medskip

The motivation of the present paper is the following.
On one hand,
the Gauss-Codazzi equations completely describe the geometry.
On the other hand,
the \PVI\ equation is complete in the classical sense \cite{PaiBSMF},
i.e.~it is impossible to add any term to \PVI\ 
under penalty of losing the Painlev\'e property
(singlevaluedness of the general solution
near all those singularities which depend on the initial conditions).
Comparing the total number of arbitrary constants in \PVI\
(six, i.e.~two movable constants plus four fixed constants $\theta_j$, see Eq.~(\ref{eqPVI}))
and the small difference (one) between the
sum of the differential order of (\ref{eqGaussCodazziR3c}) (four) 
and the number of its parameters (one, $c$),
it is natural to search for a solution in terms of the full \PVI,
i.e.~to remove the two constraints among the four $\theta_j$
in the two existing \PVI\ solutions \cite{Bonnet1867,H1897,BEK1997,BE2000}. 
\smallskip

Since surfaces described by Painlev\'e equations
ultimately arise from some ordinary differential equation (ODE),
in the present paper
we systematically look for reductions of the PDE system (\ref{eqGaussCodazziR3c})
to a system of ODEs.
\medskip

The paper is organized as follows.

In section \ref{sectionPrevious_solutions},
we recall for later reference the two solutions, found up to now, 
of the PDE system (\ref{eqGaussCodazziR3c}) in terms of \PVI.

In section \ref{sectionLie},
we establish the Lie symmetries of (\ref{eqGaussCodazziR3c})
and define essentially one reduction to a system of ODEs,
the reduced independent variable being
$\xi=G_1(z)-G_2(\barz)$, with $G_1$ and $G_2$ arbitrary.
The reduced system, which depends on one additional parameter $g$,
is further integrated with the isothermic constraint, 
i.e.~$Q$ real in some local coordinate.
The generic case $g^2\not=1$ 
provides a new solution,
in which $e^u H/Q$ and the logarithmic derivatives of $H,Q$
with respect to $\xi$
are \PVI\ functions.
The two nongeneric cases $g=-1$ and $g=1$ allow us to recover,
therefore in a systematic way, the two previous \PVI\ solutions.

\section{Previous solutions involving \PVI} 
\label{sectionPrevious_solutions}

In the point of view adopted here,
the system (\ref{eqGaussCodazziR3c}) is a set of three coupled 
partial differential equations (PDEs) in four unknowns $(u,H,Q,\barQ)$,
and, in contrast to the geometric approach,
these four unknowns will be considered as \textit{independent} and \textit{complex}.

Since the Gauss-Codazzi equations are underdetermined,
one can impose an additional constraint on the variables of (\ref{eqGaussCodazziR3c}). 
Before presenting these constraints and listing the resulting ODE systems 
defined in the complex domain,
let us recall some invariance properties and give two definitions (Bonnet surface, isothermic surface), 
or more precisely two characteristic properties.

This nonlinear system (\ref{eqGaussCodazziR3c}) admits two types of invariance.
The first one is the conformal invariance
\begin{eqnarray}
& &
\forall G_1, G_2:\ (z,\barz,e^u,H,Q,\barQ) 
\nonumber\\ & & \phantom{12345}
   \to (G_1(z),G_2(\barz),G_1'(z) G_2'(\barz) e^u, H, {G_1'(z)}^2 Q, {G_2'(\barz)}^2 \barQ). 
\label{eqLConformal}
\end{eqnarray}
The second one only exists under the condition $Q-\barQ=c$,
this is the involution 
\cite[Eq.~(4.4)]{B1994}
\cite[p.~77]{BE2000}
\cite[\S 3 p.~6]{Springborn} 
\begin{eqnarray}
& & 
(u,H,Q,\barQ) \to \left(-u,2 Q-c=2 \barQ+c,\frac{H+c}{2}, \frac{H-c}{2} \right).
\label{eqGC-Involution}
\end{eqnarray}

\begin{theorem} \cite[Theorem 2.1]{ChenLi1997} 
A Bonnet surface in $\Rcubec$ is characterized by the two conditions
\begin{eqnarray}
& & {\hskip -16.0truemm}
(\log     Q)_{z\barz}-(\log     Q)_{\barz}(\log \barQ)_{    z}=0,\ 
(\log \barQ)_{z\barz}-(\log \barQ)_{    z}(\log     Q)_{\barz}=0. 
\label{eqBonnet-surface}
\end{eqnarray}
\end{theorem}
This was proven for $c=0$ in \cite{Graustein1924} 
and for any $c$ in \cite{Tribuzy}. 

\begin{theorem} \cite[Theorem 3.2]{ChenLi1997} 
An isothermic surface in $\Rcubec$ is characterized by the condition
\begin{eqnarray}
& &
(\log Q - \log\barQ)_{z\barz}=0.
\label{eqIsothermic-surface}
\end{eqnarray}
\end{theorem}

Since the relation (\ref{eqIsothermic-surface}) is a consequence 
of the two relations (\ref{eqBonnet-surface}),
any Bonnet surface is isothermic, but not conversely,
see an example section \ref{section-iH-harmonic}.

The notation $c$ with a subscript ($\cz,\cu,\ch,\cq$, etc) 
 denotes a complex-valued constant.

\vfill \eject
\subsection{Constraint to be a Bonnet surface} 
\label{section-iQ-harmonic}

Since any Bonnet surface is isothermic,
if it is in addition umbilic-free, 
there exists an isothermic coordinate $z$ such that $1/Q(z,\barz)$ is harmonic.
Let us therefore impose the constraint that $1/Q$ be harmonic, 
\begin{eqnarray}
& & \frac{1}{Q}=f(z) + \bar f(\barz).
\label{eqQharmonic} 
\end{eqnarray}
This will lead to a solution in terms of \PVI,
first found when $c=0$ by 
Bonnet in 1867 \cite[p.~84]{Bonnet1867} 
(although often attributed to Hazzidakis \cite[p.~48]{H1897}), 
and later extrapolated to arbitrary $c$ by Bobenko and Eitner \cite[\S 4.5]{BE2000}.
\smallskip

If instead of $z$ one considers the conformal coordinate $w$ defined by
\begin{eqnarray}
& & 
\frac{\D w}{\D z}= \frac{\D z}{\D f}, 
\end{eqnarray}
then both $e^u$ and $H$ only depend on the real part of $w$,
thus defining a coupled system of two ODEs with coefficients depending on $Q$,
\begin{eqnarray}
& & {\hskip -16.0truemm}
    Q=\frac{2 \cq}{\ch} 
 \frac{4 \cz \sinh (2 \cz \bar w)}{\sinh(2 \cz      w) \sinh(4 \cz \Re(w))},\
\barQ=\frac{2 \cq}{\ch} 
 \frac{4 \cz \sinh (2 \cz      w)}{\sinh(2 \cz \bar w) \sinh(4 \cz \Re(w))},
\label{eqBE1998Qw} 
\\ & & {\hskip -16.0truemm}
e^u=4 \frac{\cq}{\ch^2} v(\Re(w)),\ H=\ch h(\Re(w)).
\label{eqBE1998uHw} 
\end{eqnarray}

Therefore the relations (\ref{eqBE1998uHw})
define a closed system of ODEs in the independent variable $\xi=\Re(w)$,
\begin{eqnarray}
& & {\hskip -14.0truemm}
\left\lbrace 
\begin{array}{ll}
\displaystyle{
(\log v)'' +8 \cq 
 \left[(h^2 -C) v 
       - \left(\frac{4 \cz}{\sinh(4 \cz \xi)}\right)^2 v^{-1}\right]=0,\
       C=\left(\frac{\cq c}{\ch}\right)^2,
}\\ \displaystyle{
v h' + \left(\frac{4 \cz}{\sinh(4 \cz \xi)}\right)^2=0,
}
\end{array}
\right.
\label{eqHazzi-Reducvh} 
\end{eqnarray}
\medskip
which admits a first integral 
(see Appendix \ref{sectionReducedGW} for methods to find first integrals)
\begin{eqnarray}
& & {\hskip -16.0truemm}
K= 
\left(\frac{v'}{v}\right)^2
 + 16 \cq 
\left[(h^2 -C) v
      + \frac{(4 \cz)^2}{\sinh^2(4 \cz \xi)} v^{-1}
      -2 (4 \cz \coth(4 \cz \xi)) h\right].
\label{eqHazzi-firsthv} 
\end{eqnarray}
The elimination of $v$ from (\ref{eqHazzi-Reducvh})${}_2$
yields a third order first degree ODE
\cite[\S 11 p 84 Eq (52)]{Bonnet1867}, 
\begin{eqnarray}
& &
(\log h')''
+ 8 \cq h'
- \left(\frac{4 \cz}{\sinh(4 \cz \xi)}\right)^2
 \left(8 \cq \frac{h^2-C}{h'} +2 \right)=0,
\label{eqBonnet-ODE3h}
\end{eqnarray} 
and its first integral \cite[p.~48]{H1897} \cite[Eq.~(4.20)]{BE2000} 
\begin{eqnarray}
& & {\hskip -15.0 truemm}
K= 
\left(\frac{h''}{h'} +8 \cz \coth(4 \cz \xi)\right)^2
\nonumber \\ & & {\hskip -15.0 truemm} \phantom{1234}
+16 \cq \left[
             \left(\frac{4 \cz}{\sinh(4 \cz \xi)}\right)^2 \frac{h^2 -C}{h'}
             + h'
             + 8 \cz \coth(4 \cz \xi) h\right].
\label{eqHazzi-ODE2h} 
\end{eqnarray}

Therefore, the six arbitrary constants of Bonnet surfaces 
\cite[p.~85]{Bonnet1867} 
\cite[pp.~55]{Cartan1942} 
\cite{Chern1985}
are the two origins of $z$, $\barz$, the two first integrals $\cz$, $K$,
and the two constants of integration of (\ref{eqHazzi-ODE2h}). 

The integration of this second order second degree ODE
could of course not be performed by Bonnet nor Hazzidakis 
since the complete \PVI\ would only be discovered by R.~Fuchs in 1905. 
This integration, first achieved in \cite{BE1998}, is as follows.

When $\cz\not=0$,
the ODEs (\ref{eqHazzi-ODE2h}) and SD-Ia (\ref{eqSDIa}) 
obeyed by anyone of the 24 Hamiltonians of \PVI\
(Appendix \ref{AppendixPVI})
are equivalent under a homographic transformation,
e.g.
\begin{eqnarray}
& & {\hskip -15.0 truemm}
\left\lbrace 
\begin{array}{ll}
\displaystyle{
h=\frac{2 \cz}{\cq} X(X-1) \Hu_{\rm VI},\  
X=\frac{1}{1-e^{8 \cz \xi}},\
}\\ \displaystyle{
A_0=\frac{K}{(8 \cz)^2},
A_2=A_4=0,\ 
A_3=\frac{\cq^2 C}{\cz^2}-\frac{K^2}{4(8 \cz)^2},
}
\end{array}
\right.
\label{eqHazzi-IntegrationP6} 
\end{eqnarray}
in which $\Hu_{\rm VI}$ is defined in Eq.~(\ref{eqHamVICu}).
The monodromy exponents of \PVI\ are then, for instance (see (\ref{SDIanjAjtj})),
\begin{eqnarray}
& & {\hskip -15.0 truemm}
\left\lbrace 
\begin{array}{ll}
\displaystyle{
(\theta_\infty^2,\theta_0^2,\theta_1^2,\Theta_X^2)
=(0,t_m+\sqrt{t_m^2+t_u},t_m-\sqrt{t_m^2+t_u},0),\
}\\ \displaystyle{
K=(8 \cz)^2 t_m,\
C=\frac{\cz^2}{4 \cq^2} (t_m^2+t_u). 
}
\end{array}
\right.
\label{eqHazzitheta} 
\end{eqnarray}

When $\cz=0$, 
Eq.~(\ref{eqHazzi-ODE2h}) becomes
\begin{eqnarray}
& & {\hskip -15.0 truemm}
K= 
\left(\frac{h''}{h'}+\frac{2}{\xi}\right)^2
+16 \cq \left[\frac{h^2 -C}{\xi^2 h'} + h' + 2 \frac{h}{\xi}\right],
\label{eqHazzi-ODE2h0} 
\end{eqnarray}
and $h$ is then an affine function of the Hamiltonian of \PV\ \cite{CosScou},
\begin{eqnarray}
& & {\hskip -16.0truemm}
-\left(X Y''\right)^2 - 4 {Y'}^2 (X Y'-Y)
+ A_1 (X Y'-Y)^2
+ A_2 (X Y'-Y)
+ A_3 Y'
+ A_4,
\nonumber\\ & & {\hskip -16.0truemm}
h=\frac{Y_{\rm V}}{4 \cq \xi},\ \xi=X,\
A_1=K,\ A_2=-(8 \cq )^2 C,\ A_3=0,\ A_4=0.
\label{eqHazzi-IntegrationP5} 
\end{eqnarray}

\vfill\eject

\subsection{Constraint to be a surface dual to a Bonnet surface} 
\label{section-iH-harmonic}

The involution (\ref{eqGC-Involution})
maps an umbilic-free Bonnet surface 
(i.e.~there exists an isothermic coordinate $z$ such that $1/Q(z,\barz)$ is harmonic) 
to another surface which is not a Bonnet surface but 
which is dual to a Bonnet surface \cite[Prop.~4.7.1 page 77]{BE2000}
and has a harmonic inverse mean curvature $1/H$ (HIMC).
\smallskip

Let us therefore impose the constraint that $1/H$ be harmonic, 
\begin{eqnarray}
& &
(1/H)_{z \barz}=0. 
\label{eqHharmonic} 
\end{eqnarray}
Initially presented in \cite[section 7]{B1994},
this constrained system has been integrated with \PVI\ in \cite[section 5]{BEK1997}. 
After some algebra, 
one obtains
\begin{eqnarray}
\frac{1}{H}
=\frac{2 \cz}{\ch}\left(\coth(2 \cz z)+\coth(2 \cz \barz)\right)
=\frac{2 \cz}{\ch}\frac{\sinh(4 \cz x)}
                       {\sinh(2 \cz z) \sinh(2 \cz \barz)},                     
\label{eqBEK-H} 
\end{eqnarray}
$c$ must be equal to zero,
then $e^u H^2$ and $Q$ are proven to only depend on $x$ \cite{BEK1997},
\begin{eqnarray}
& & {\hskip -15.0 truemm}
\left\lbrace 
\begin{array}{ll}
\displaystyle{
 e^u=\frac{4 \cu v(x)}{H^2},
}\\ \displaystyle{
    Q=\frac{4\cq}{\ch}\left(\frac{2\cz}{\sinh(2\cz x)}\right)^2 (q(x)+i\theta),\
\barQ=\frac{4\cq}{\ch}\left(\frac{2\cz}{\sinh(2\cz x)}\right)^2 (r(x)-i\theta),
}
\end{array}
\right.
\label{eqBEK-uQ} 
\end{eqnarray} 
in which $\theta$ is an arbitrary constant.
The reduced system \cite{BEK1997}, 
\begin{eqnarray}
& & {\hskip -15.0 truemm}
\left\lbrace 
\begin{array}{ll}
\displaystyle{
(\log v)''
+\frac{8}{\cu} v
-\left(\frac{4 \cz}{\sinh(4 \cz x)}\right)^2
 \left(8 \cq^2\frac{(q+i \theta)(r-i \theta)}{\cu v} +2 \right)=0,\
}\\ \displaystyle{
\cu v-\cq q'=0,
}\\ \displaystyle{
\cu v-\cq r'=0,
}
\end{array}
\right.
\label{eqBEK-Reducqrv}
\end{eqnarray}
admits the first integral
\begin{eqnarray}
& & {\hskip -23.0 truemm}
K= 
\left(\frac{v'}{v} +8 \cz \coth(4 \cz x)\right)^2
\nonumber\\ & & {\hskip -23.0 truemm} \phantom{12}
+16 \left(\frac{4 \cz}{\sinh(4 \cz x)}\right)^2 
    \frac{\cq^2}{\cu} (q+i \theta)(r-i \theta)  
		+16 \cu v
		+16 \cq \frac{4 \cz}{\sinh(4 \cz x)} (q + r).
\label{eqBEK-firstqrv} 
\end{eqnarray}

Without loss of generality, one can choose zero for the constant value of $q-r$.
The elimination of $v$ from (\ref{eqBEK-Reducqrv})${}_2$
yields a third order first degree ODE for $q(x)$,
\begin{eqnarray}
& & {\hskip -13.0truemm}
(\log q')''
+ 8 \cq q'
-\left(\frac{4 \cz}{\sinh(4 \cz x)}\right)^2
 \left(8 \cq \frac{q^2 + \theta^2}{q'} +2 \right)=0,\
\label{eqBEK-ODE3q} 
\end{eqnarray} 
which admits the first integral
\begin{eqnarray}
& & {\hskip -10.0 truemm}                                  
K= 
\left(\frac{q''}{q'} +8 \cz \coth(4 \cz x)\right)^2
\nonumber\\ & &  {\hskip -10.0 truemm} \phantom{12}
+16 \cq \left[
    \left(\frac{4 \cz}{\sinh(4 \cz x)}\right)^2 \frac{q^2+\theta^2}{q'}
    + q'
    + 8 \cz \coth(4 \cz x) q\right].
\label{eqBEK-ODE2q} 
\end{eqnarray}
Due to the above-mentioned involution (\ref{eqGC-Involution}),
the second order second degree ODEs
(\ref{eqHazzi-ODE2h}) and (\ref{eqBEK-ODE2q})
are exchanged under 
\begin{eqnarray}
& &
h \to q,\
C \to - \theta^2.
\end{eqnarray}
Therefore (\ref{eqBEK-ODE2q}) integrates with either \PVI\ or \PV\,
depending on whether $\cz$ is nonzero or zero,
see previous section \ref{section-iQ-harmonic}.

\vfill\eject

\section{Lie point symmetries and reductions} 
\label{sectionLie}

By requiring the Gauss-Codazzi equations to be invariant under a Lie point symmetry,
one can define a reduction to a system (again underdetermined)
of three coupled ODEs in four dependent variables.
Then, as mentioned in section \ref{sectionPrevious_solutions},
an additional constraint allows us to integrate this system in terms of \PVI\ functions. 

\subsection{Lie point symmetries of the Gauss-Codazzi system}

Lie point symmetries are obtained by a classical computation presented for instance in
\cite{OlverBook,OvsiannikovBook}.

Denoting for convenience $U=e^u, R=\barQ$,
the result for the system (\ref{eqGaussCodazziR3c}) is the 
infinite-dimensional Lie algebra generated by
\begin{eqnarray}
& &
{\hskip -15.0 truemm}
\left\lbrace 
\begin{array}{ll}
\displaystyle{
X(F)=F(    z) \partial_z     + F'(    z) (-2 Q \partial_Q-U \partial_U),
}\\ \displaystyle{
Y(G)=G(\barz) \partial_\barz + G'(\barz) (-2 R \partial_R-U \partial_U),
}\\ \displaystyle{
(c=0 \hbox{ only})\ a=-H \partial_H+ Q \partial_Q+ R \partial_R+2 U \partial_U,\
}
\end{array}
\right.
\label{eqGC_Lie-algebra-infinite} 
\end{eqnarray} 
in which $F,G$ are arbitrary functions of one variable.
Its table of commutation is
\begin{eqnarray}
& &
{\hskip -15.0 truemm}
\left\lbrace 
\begin{array}{ll}
\displaystyle{
[X(F_1),X(F_2)]=X(F_1 F_2'-F_1' F_2),\
}\\ \displaystyle{
[Y(G_1),Y(G_2)]=Y(G_1 G_2'-G_1' G_2),\
}\\ \displaystyle{
[X(F),Y(G)]=0,\
[X(F),a]=0,\
[Y(G),a]=0.
}
\end{array}
\right.
\label{eqGC_Lie-algebra-infinite-commu} 
\end{eqnarray} 
The largest finite-dimensional subalgebra is
defined by the seven generators $X(z^j),Y(\barz^j)$, $j=0,1,2$ and $a$,
\begin{eqnarray}
& &
{\hskip -15.0 truemm}
\left\lbrace 
\begin{array}{ll}
\displaystyle{
e_0=    \partial_z,\
e_1=z   \partial_z        -2 Q \partial_Q-U \partial_U,\
e_2=z^2 \partial_z + 2 z (-2 Q \partial_Q-U \partial_U),
}\\ \displaystyle{
f_0=    \partial_\barz,\
f_1=\barz   \partial_\barz            -2 R \partial_R-U \partial_U,\
f_2=\barz^2 \partial_\barz + 2 \barz (-2 R \partial_R-U \partial_U),
}\\ \displaystyle{
(c=0 \hbox{ only})\ a=-H \partial_H + Q \partial_Q + R \partial_R + 2 U \partial_U,\
}
\end{array}
\right.
\label{eqGC_Lie-algebra-finite} 
\end{eqnarray} 
with the nonzero commutators
\begin{eqnarray}
& &
[e_0,e_1]=e_0,\ [e_0,e_2]=2 e_1,\ [e_1,e_2]=e_2,\ 
\nonumber \\ & &
[f_0,f_1]=f_0,\ [f_0,f_2]=2 f_1,\ [f_1,f_2]=f_2,\ 
\label{eqGaussCodazziCommuR3} 
\end{eqnarray}
This seven-dimensional algebra is therefore the direct sum
\begin{eqnarray}
& &
\lbrace{e_0,e_1,e_2\rbrace} \directsum \lbrace{f_0,f_1,f_2\rbrace} \directsum \lbrace{a \rbrace},
\end{eqnarray}
and the subalgebra $\lbrace{e_0,e_1,e_2\rbrace}$ is isomorphic to su(1,1),
see \cite[p.~1451]{PW}. 

There is essentially one resulting reduction,
in which the reduced variable is the sum (or, equivalently, the product)
of an arbitrary function of $z$ and an arbitrary function of $\barz$.

\vfill\eject
\subsection{Reduction} 
\label{sectionReduc-Lwhat}

The infinite-dimensional Lie algebra (\ref{eqGC_Lie-algebra-infinite})
defines the reduction
(with the notation $g_1=F, g_2=G$),
in which $a_1,a_2$ denote arbitrary complex constants,
\begin{eqnarray}
& & {\hskip -15.0 truemm}
\left\lbrace 
\begin{array}{ll}
\displaystyle{
\xi=\log g_1(z) - \log g_2(\barz),\
e^u=g_1^{2 a_1-1} g_2^{2 a_2-1} g_1' g_2' \tildev,\
H  =g_1^{- a_1}   g_2^{- a_2}             \tildeh,\
}\\ \displaystyle{
Q  =g_1^{  a_1-2} g_2^{  a_2} {g_1'}^2    \tildeq,\
R  =g_2^{  a_2-2} g_1^{  a_1} {g_2'}^2    \tilder,\
}
\end{array}
\right.
\label{eqGCReduca1a2}
\end{eqnarray}
which is indeed a reduction iff the product $\as c$ is zero,
with the notation $\as=a_1+a_2, \ad=a_1-a_2$.
The reduced system in $(\tildev,\tildeh,\tildeq,\tilder)(\xi)$,
\begin{eqnarray}
& & 
\as c=0:\
\left\lbrace 
\begin{array}{ll}
\displaystyle{
(\log \tildev)'' - \frac{1}{2}\left(\tildeh^2-c^2 e^{\displaystyle{\ad \xi}}\right) \tildev 
 + 2 \tildeq \tilder \tildev^{-1}=0,\
}\\ \displaystyle{
2 \tildeq' + \tildev \tildeh' - \frac{1}{2} (\as+\ad) \tildev \tildeh -(\as-\ad) q=0,
}\\ \displaystyle{
2 \tilder' + \tildev \tildeh' + \frac{1}{2} (\as-\ad) \tildev \tildeh +(\as+\ad) r=0,
}
\end{array}
\right.
\label{eqGCsystemra1a2}
\end{eqnarray} 
whose reduced moving frame equations 
are given in the Appendix section \ref{sectionReducedGW},
admits the first integral
\begin{eqnarray}
& & {\hskip -23.0 truemm}
K 
= \left(\frac{\tildev'}{\tildev}+\ad\right)^2
-  \left(h^2-c^2 e^{\displaystyle{\ad \xi  }} \right) \tildev 
-2 \tildeh (\tildeq+\tilder) - 2 c e^{\displaystyle{\ad \xi/2}} (\tildeq-\tilder) 
- 4 \frac{\tildeq \tilder}{\tildev}\cdot
\label{eqGCfirsta1a2} 
\end{eqnarray} 
Except when $\ad c \not=0$, 
the system (\ref{eqGCsystemra1a2}) and its first integral
are autonomous.

Under the change of variables defined by
\begin{eqnarray}
& &  {\hskip -20.0 truemm}
e^\xi=\eta,\
\tildev=e^{(2 a_2-1)\xi } v,\
\tildeh=e^{ - a_2   \xi } h,\
\tildeq=e^{   a_2   \xi } q,\
\tilder=e^{(  a_2-2)\xi } r,\
\as=g+1,
\label{eqGC_systemg-xi_to_eta}
\end{eqnarray} 
the system
(\ref{eqGCsystemra1a2})--(\ref{eqGCfirsta1a2}) 
is mapped to a system with one less parameter,
\begin{eqnarray}
& &
(1+g) c=0:\ 
\left\lbrace 
\begin{array}{ll}
\displaystyle{
\left(\eta \frac{v'}{v} \right)'-\frac{h^2-c^2}{2\eta^2}v+2\frac{q r}{\eta^2 v}=0,
}\\ \displaystyle{
-\eta v h' -2 \eta^2 q' +(1+g) v h=0,
}\\ \displaystyle{
-\eta^2 v h' -2 \eta r' +2 (1-g) r=0,
}
\end{array}
\right.
\label{eqGC_systemg} 
\end{eqnarray} 
which admits the first integral
\begin{eqnarray}
& & {\hskip -15.0 truemm}
K=
\left(\eta \frac{v'}{v} +g \right)^2
-\frac{(h^2-c^2) v}{\eta}
-2 (h+c) q
-2 \frac{(h-c)}{\eta^2} r
-4 \frac{q r}{\eta v}\cdot
\label{eqGC_firstg} 
\end{eqnarray} 

{}From now on, we impose the constraint $q=r$ 
(equivalent to the isothermic constraint $Q=\barQ=R$ 
under a suitable choice of the so far arbitrary functions $g_1(z), g_2(\barz)$).

\vfill\eject
\subsubsection{Reduction, generic isothermic case, $g^2\not=1, c=0$ } 
\label{sectionReduc-Lwhat-generic}

In the generic case $g^2\not=1$, $c$ must vanish and
the three fields $h v/q,h'/h,q'/q$ 
are equivalent under the homographic group,
\begin{eqnarray}
& & {\hskip -17.0 truemm}
g^2\not=1:\ 
1+\frac{\eta h v}{2 q}
=\frac{2 - (1-\eta^2) \eta h'/h}{(1+g)}
=1-\frac{(1-g) \eta^2}{(1-\eta^2) \eta q'/q-(1-g)},
\label{eqHomo3}
\end{eqnarray} 
therefore it is more convenient to choose them 
as dependent variables rather than $v,h,q$.
After elimination of two of these three variables $h v/q,h'/h,q'/q$,
the remaining variable is found to obey a second order first degree ODE.
Since the three such ODEs only differ by the homographic transformation (\ref{eqHomo3}), 
it is sufficient to consider only one of them. 
In the case $W=\eta h v/(2 q)$ for instance, this is, 
\begin{eqnarray} 
& & {\hskip -15.0 truemm}
\frac{\D^2 W}{\D \eta^2}=
 \frac{1}{2} \left[\frac{1}{W} + \frac{1}{W+1} + \frac{1}{W+\eta^2} \right] \left(\frac{\D W}{\D \eta}\right)^2
\nonumber\\ & & {\hskip -15.0 truemm}
- \left[\frac{1}{\eta} - \frac{2 \eta}{1-\eta^2} - \frac{2 \eta}{W+\eta^2} \right] \frac{\D W}{\D \eta}
+ \frac{2 W (W+1) (W+\eta^2)}{\eta^2 (1-\eta^2)^2} \times
\nonumber\\ & & {\hskip -15.0 truemm}
  \left[
	  \left(\frac{1+g}{2} \right)^2 
	- \left(\frac{1-g}{2} \right)^2 \frac{\eta^2}{4 W^2} 
	+ \frac{K}{4} \frac{\eta^2-1}{(W+1)^2}
  + (1-\frac{K}{4}) \frac{\eta^2 (\eta^2-1)}{(W+\eta^2)^2} \right].
\label{eqfvhq}
\end{eqnarray} 
This ODE belongs to the class studied by Painlev\'e and Gambier \cite{GambierThese}
\begin{eqnarray}
& &
W''=A_2(W,\eta) {W'}^2+A_1(W,\eta) W' +A_0(W,\eta),
\label{eqloghprime}
\end{eqnarray} 
the rational function $A_2$ of $W$ is the sum of four simple poles 
with equal residues $1/2$,
located at
\begin{eqnarray}
& & W=\infty,\ 0,\ -1,\ -\eta^2, 
\end{eqnarray} 
these four poles are distinct and their crossratio is not a constant,
therefore the ODE (\ref{eqfvhq}) is a homographic transform of \PVI, 
Eq.~(\ref{eqPVI}).
This provides the general solution in terms of a \PVI\ function $V(X)$ 
with two constraints among its four parameters $\theta_j$,
i.e.~as many as in the two existing solutions,
\begin{eqnarray}
& & {\hskip -11.0 truemm}
\left\lbrace 
\begin{array}{ll}
\displaystyle{
W=\frac{\eta h v}{2 q}=\frac{V-X}{X},\ 
}\\ \displaystyle{
\frac{\D \log h}{\D X}=\frac{1+g}{2(X-1)} - \frac{1-g}{2(V-X)},\  
\frac{\D \log q}{\D X}=\frac{1  }{  X-1 } - \frac{(1+g) V}{2 X(X-1)},\  
}\\ \displaystyle{
h q=\frac{2(X-1)}{V(V-1)(V-X)} 
}\\ \displaystyle{ \phantom{12345} \times
\left[\left(X(X-1) V' + \frac{1-g}{2}V(V-1) \right)^2 -\frac{K}{4} (V-X)^2 \right], 
}\\ \displaystyle{
X=\frac{1}{1-\eta^2},\
(\theta_\infty^2,\theta_0^2,\theta_1^2,\theta_X^2)=
\left(\left(\frac{1+g}{2} \right)^2, 
\frac{K}{4}, \frac{K}{4}, \left(\frac{1-g}{2} \right)^2\right).
}
\end{array}
\right.
\label{eqGC_Lie_red_g_solutionP6} 
\end{eqnarray} 
Several remarks are in order.
\begin{enumerate}

\item
This solution is invariant under the involution (\ref{eqGC-Involution}),
which acts as a homography on $V$.

\item
For generic values of the $\theta_j$'s,
there exists no rational function of $V,V'$
which would represent $h$ or $q$.
Only their product $h q$ is rational.

\item
Since there exists a choice of square roots 
allowing the monodromy exponents $\theta_j$ to have an integral sum,
there always exists a one-parameter solution of this \PVI\ 
which is an algebraic transform of the hypergeometric function 
and obeys, for instance,
\begin{eqnarray}
& & {\hskip -11.0 truemm}
g^2\not=1:\ 
\left\lbrace 
\begin{array}{ll}
\displaystyle{
\frac{X (X-1) V'}{V (V-1)(V-x)}
 +\frac{\theta_0}{V}+\frac{\theta_1}{V-1}+\frac{\theta_X-1}{V-X}=0,
}\\ \displaystyle{
\theta_\infty=\frac{1+g}{2},\ 
\theta_0=\frac{\sqrt{K}}{2},\ 
\theta_1=-\frac{\sqrt{K}}{2},\ 
\theta_X=\frac{1-g}{2}\cdot
}
\end{array}
\right.
\label{eqP6hyperg} 
\end{eqnarray} 

\item
The values of the monodromy exponents even allow 
algebraic solutions to exist.
For instance, the cube solution of \PVI\ \cite{KitaevP6cube}
\begin{eqnarray}
& & {\hskip -11.0 truemm}
\left\lbrace 
\begin{array}{ll}
\displaystyle{
-X^2 + 3 X V -3 X V^2 + 2 X V^3 - V^3=0
}\\ \displaystyle{
(\theta_\infty,\theta_0,\theta_1,\theta_X)=(1/3,a,a,2 a),\ a \hbox{ arbitrary},
}
\end{array}
\right.
\label{eqP6cube} 
\end{eqnarray} 
matches (\ref{eqGC_Lie_red_g_solutionP6})
for $a=2/3$, $g=-1/3$, $K=16/9$.

\item
{}From the definitions (\ref{eqGCReduca1a2}) 
and (\ref{eqGC_systemg-xi_to_eta}),
we could not build a function of $(e^u,H,Q)$
which would be a harmonic function of $(z,\barz)$,
like for instance $z^2 (\eta^2-1) / \eta^2$.
              
\item 
One may wonder whether the invariances of \PVI\ allow
mapping the present surface to a Bonnet surface. 
The answer is negative, at least for generic values of the monodromy exponents.
Indeed, the elementary birational transformation which conserves \PVI\ \cite{Okamoto1987I} 
acts as an affine transformation on the $\theta_j$'s
(see, e.g.~\cite[\S B.3.1]{CMBook})
and, starting from the $\theta_j$'s in (\ref{eqGC_Lie_red_g_solutionP6}),
one can easily make one transformed $\theta_j$ equal to zero,
but not two. 

\end{enumerate}

\vfill\eject

\subsubsection{Reduction, isothermic case, $g=1,c=0$}  
\label{sectionReduc-Lwhat-gplus1}

The elimination of $v$ among the last two equations of (\ref{eqGC_systemg})
defines two subcases,
\begin{eqnarray}
& &
q' \left[(1-\eta^2) \eta h'-2 h\right]=0.
\end{eqnarray} 

The first subcase $q'=0$ is equivalent to
\begin{eqnarray}
& &
g=1,\ c=0,\
h=0,\
q'=0,\
\left(\eta \frac{v'}{v}+1\right)^2 -K - 4 \frac{q^2}{\eta v}=0, 
\end{eqnarray} 
and its general solution
\begin{eqnarray}
& & {\hskip -18.0 truemm}
g=1,\ c=0,\
h=0,\
q=\cq,\
\eta=e^\xi,\
\eta v=\left\lbrace \begin{array}{ll} \displaystyle{
\frac{4 \cq^2}{K} \sinh^2(\frac{\sqrt{K}}{2}(\xi-\xi_0)),\ \cq K\not=0,
}\\ \displaystyle{
\cq^2 (\xi-\xi_0)^2,\ \cq\not=0,\ K=0,
}\\ \displaystyle{
v_0 e^{\sqrt{K}(\xi-\xi_0)},\ \cq=0,\
}\end{array}\right.
\end{eqnarray} 
depends on at most three movable constants among $K,\cq,\xi_0,v_0$. 

The second subcase $q'\not=0$ is equivalent to
\begin{eqnarray}
& &
{\hskip -10.0 truemm}
g=1,\ c=0,\
h=\ch \frac{\eta^2}{\eta^2-1},\
v=-\frac{q'}{h},\
\ch\not=0,
\end{eqnarray} 
in which $q$ is the Hamiltonian of \PVI, see Eqs.~(\ref{eqHPVInorm})--(\ref{eqODEHPVI}),
\begin{eqnarray}
& &
{\hskip -25.0 truemm}
q(\eta)= - \frac{8}{\ch} X(X-1) \Hu_{\rm VI},\ 
X=\frac{1}{1-\eta^2},\
A_0=\frac{K}{4},\ A_2=0,\  A_3=0,\  A_4=0,
\label{eqgp1solq}
\end{eqnarray} 
and the four integration constants are $K,\ch$ and the two movable constants of \PVI.
This solution,
for which $1/H$ is harmonic,
\begin{eqnarray}
& &
\frac{1}{H}=\frac{1}{\ch}\left(g_1^2(z)-g_2^2(\barz)\right),
\end{eqnarray} 
is the particular case $\theta=0$ 
of the HIMC solution of section \ref{section-iH-harmonic}.

\subsubsection{Reduction, isothermic case, $g=-1,c \hbox{ arbitrary}$} 
\label{sectionReduc-Lwhat-gminus1}

The elimination of $v$ among the last two equations of (\ref{eqGC_systemg})
defines two subcases,
\begin{eqnarray}
& &
h' \left[(1-\eta^2) \eta q'-2 q\right]=0.
\end{eqnarray} 

The first subcase $h'=0$ is equivalent to
\begin{eqnarray}
& & {\hskip -19.0 truemm}
g=-1,\ c \hbox{ arbitrary},\
h'=0,\
q=0,\
\left(\eta \frac{v'}{v}-1\right)^2 -K -(h^2-c^2) \frac{v}{\eta}=0, 
\end{eqnarray} 
and its general solution
\begin{eqnarray}
& & {\hskip -18.0 truemm}
g=-1,\ 
h=h_0,\
q=0,\
\frac{\eta}{v}=\left\lbrace \begin{array}{ll} \displaystyle{
\frac{h_0^2-c^2}{K} \sinh^2((\sqrt{K}/2) \xi-\xi_0),\ (h_0^2-c^2) K\not=0,
}\\ \displaystyle{
\frac{h_0^2-c^2}{4} (\xi-\xi_0)^2,\ (h_0^2-c^2)\not=0,\ K=0,
}\\ \displaystyle{
v_0 e^{\sqrt{K}(\xi-\xi_0)},\ h_0^2-c^2=0,\
}\end{array}\right.
\end{eqnarray} 
depends on at most three movable constants among $K,h_0,\xi_0,v_0$. 

The second subcase $h'\not=0$ is equivalent to
\begin{eqnarray}
& &
{\hskip -10.0 truemm}
g=-1,\ c \hbox{ arbitrary},\
q=\cq \frac{\eta^2}{\eta^2-1},\
v=\frac{4 \cq \eta^2}{(1-\eta^2)^2 h'},\
\cq\not=0,
\end{eqnarray} 
in which $h$ is an affine function of the Hamiltonian of \PVI, 
see Eqs.~(\ref{eqHPVInorm})--(\ref{eqODEHPVI}),
\begin{eqnarray}
& &
\left\lbrace \begin{array}{ll} \displaystyle{
h= -\frac{8}{\cq} X(X-1) \Hu_{\rm VI},\
X=\frac{1}{1-\eta^2},\
}\\ \displaystyle{ 
A_0=\frac{K + 2 \cq c}{4},\
A_2=0,\
A_3=-\frac{K(K+4 \cq c)}{64}
A_4=0,
}\end{array}\right.
\label{eqgm1solh}
\end{eqnarray} 
and the four integration constants are $K,\cq$ and the two movable constants of \PVI.
This solution, in which $1/Q$ and $1/R$ are harmonic up to some conformal transformation,
\begin{eqnarray}
& &
\frac{\left(g_1'(    z)\right)^2}{Q}=
\frac{\left(g_2'(\barz)\right)^2}{R}=\frac{1}{\cq}\left(g_1^2(z)-g_2^2(\barz)\right),\
\end{eqnarray} 
is identical to the solution of Bonnet described in section \ref{section-iQ-harmonic}. 
\medskip

These two nongeneric solutions $g=1,-1$ are of course exchanged 
under the involution (\ref{eqGC-Involution}).
\medskip
 
In Table \ref{TableComplex},
we have gathered all the (complex) solutions of 
the $\Rcubec$ Gauss-Codazzi equations
which involve \PVI\ or \PV.
Let us remind once again that all constants $\cz, \cu, \cq,$ etc,
are complex.
For instance, $4 \cz/\sinh(4 \cz x)$ denotes at the same time
$1/\sinh x$, $1/\sin x$ or $1/x$.

\vfill\eject

\tabcolsep=1.5truemm
\tabcolsep=0.5truemm

\vspace{20pt}
\begin{landscape}
\begin{table}[h] 
\caption[Complex solutions of the $\Rcubec$ Gauss-Codazzi equations.]{
         Complex solutions in terms of Painlev\'e \PVI\ or \PV\ functions. \hfill\break
This table collects both old and new such solutions.
Other solutions (in terms of classical functions) can be found in the text.
The second column displays the function 
which is harmonic.
The column ``Link'' refers to the explicit link to a Painlev\'e function given in the text.
The column ``arb'' lists the arbitrary constants in \PVI\ or \PV.
}
\vspace{0.2truecm}
\begin{center}
\begin{tabular}{| r | l | l | l | l | l | l | l | l | l |}
\hline 
Reduction & Harmo & $e^u$ & $H$ & $Q$ & Integration & $c$ & arb & Link & Ref \\ \hline \hline
$x$ & $1/Q$
    & $\left(\frac{4 \cz}{\sinh(4 \cz x)}\right)^2 {h'(x)}^{-1}$
    & $h(x)$ 
    & $\frac{4 \cz \sinh (2 \cz \barz)}{\sinh(2 \cz     z) \sinh(4 \cz x)}$ 
    & $h(x)=\Hu_{\rm VI}$ 
		& $c$ 
    & $K,c$
    & (\ref{eqHazzi-IntegrationP6}) 
    & \cite{Bonnet1867} \\ \hline 
$x$ & $1/Q$ 
    & $x^{-2} {h'(x)}^{-1}$
    & $h(x)$ 
    & $ \barz/(z x)$ 
    & $h(x)=\Hu_{\rm V}$ 
		& $c$
    & $K,c$
    & (\ref{eqHazzi-IntegrationP5})
		& \cite{Bonnet1867} \\ \hline 
    \hline
$x$ & $1/H$ 
    & $q'(x) H^{-2}$
    & $\ch \frac {\sinh (2 \cz \barz) \sinh(2 \cz z)}
                 {2 \cz \sinh(4 \cz x)}$ 
    & $\left(\frac{2\cz}{\sinh(2\cz x)}\right)^2 (q+i\theta)$
    & $q(x)=\Hu_{\rm VI}$ 
		& $0$
    & $K,\theta$
    & (\ref{eqHazzi-IntegrationP6}) 
    & \cite{BEK1997} \\ \hline 
$x$ & $1/H$
    & $q'(x) H^{-2}$
    & $z \barz/(2 x)$ 
    & $x^{-2} (q(x)+i\theta)$
    & $q(x)=\Hu_{\rm V}$ 
		& $0$ 
    & $K,\theta$
    & (\ref{eqHazzi-IntegrationP5})
    & \cite{BEK1997} \\ \hline 
    \hline
$\eta=z/\barz$ 
    & ?
    & $z^{2 g}  v(\eta)$
    & $z^{-1-g} h(\eta)$
    & $z^{-1+g} q(\eta)$
    & Reduction (\ref{eqGC_systemg})
		& $0$ 
    &
    &
    &            \\ \hline 
$g=1$
    & $1/H$
    & $-z^{2} q'(\eta) \frac{\eta^2-1}{\eta^2}$
    & $z^{-2} \eta^2/(\eta^2-1)$
    & $q(\eta)$
    & $q(\eta)=\Hu_{\rm VI}$ 
	  & $0$ 
    & $K$
    & (\ref{eqgp1solq})
    & \cite{BEK1997} \\ \hline 
$g=-1$
    & $1/Q$
    & $z^{-2} \frac{\eta^2}{(1-\eta^2)^2 h'(\eta)}$
    & $        h(\eta)$
    & $z^{-2} \frac{\eta^2}{\eta^2-1}$ 
    & $h(\eta)=\Hu_{\rm VI}$ 
	  & $c$ 
    & $K,c$
    & (\ref{eqgm1solh})
    & \cite{Bonnet1867} \\ \hline 
$g^2\not=1$
    & ?
    & $z^{2 g}  v(\eta)$
    & $z^{-1-g} h(\eta)$
    & $z^{-1+g} q(\eta)$
    & $h'/h,q'/q, h v/q =\PVI$ 
		& $0$ 
    & $K,g$
    & (\ref{eqGC_Lie_red_g_solutionP6})
    & this paper \\ \hline 
\end{tabular}
\end{center}
\label{TableComplex}
\end{table}

\end{landscape}

\vfill\eject

\subsection{On the linear representations of \PVI} 
\label{sectionGWPVI}

In each of the three cases $g^2\not=1$, $g=1$, $g=-1$, 
the reduced moving frame equations (\ref{eqGW-xi-t}) 
define a linear representation of a codimension-two \PVI.
For $g^2=1$, these have been shown \cite[\S 3.4.1 and 4.9.1]{BE2000}
to be essentially not different from the isomonodromic matrix Lax pairs
of Jimbo and Miwa,
the spectral parameter $t$ being defined by (\ref{eqNOTspectralparam}).

In the case $g^2\not=1$, 
the reduced moving frame equations (\ref{eqGW-xi-t})
do not depend on the parameter $t$
and therefore define a linear representation of \PVI\
which is not a Lax pair.
This traceless, well-balanced, representation is,
\begin{eqnarray}
& & {\hskip -12.0truemm}
\D \Psi=\tilde U_r \Psi \D \xi + \tilde V_r \Psi \D t
       =       U_r \Psi \D  X  +        V_r \Psi \D t,
\nonumber\\ & & {\hskip -12.0truemm}
U_r=\frac{1}{4 X(X-1)}\left[\frac{1+g}{2}\pmatrix{1 & 0 \cr 0 & -1 \cr}\right.
\nonumber\\ & & 
\left.+ \frac{\sqrt{R_{+} R_{-}}}{\sqrt{V(V-1)}(V-X)}  
\pmatrix{0 & V-2 X+1 \cr V-2 X & 0 \cr}
\right],
\nonumber\\ & & {\hskip -12.0truemm}
V_r=\frac{X(X-1)V'-V(V-1)(1+g)/2}{2(V-X)} \pmatrix{1 & 0 \cr 0 & -1 \cr}
\nonumber\\ & & 
+\frac{\sqrt{R_{+} R_{-}}}{2 \sqrt{V(V-1)}(V-X)}  
\pmatrix{0 & -(V-1) \cr V & 0 \cr},
\nonumber\\ & & {\hskip -12.0truemm}
R_{\pm}=X(X-1) \frac{\D V}{\D X} + V(V-1)(V-X)
\left[\pm \frac{\theta_{0}}{V} \mp \frac{\theta_{1}}{V-1}+\frac{\theta_{X}-1}{V-X} \right],
\nonumber\\ & & {\hskip -12.0truemm}
(\theta_\infty^2,\theta_0^2,\theta_1^2,\theta_X^2)=
\left(\left(\frac{1+g}{2} \right)^2, 
\frac{K}{4}, \frac{K}{4}, \left(\frac{1-g}{2} \right)^2\right),
\label{eqGW-X-t-codim2}
\end{eqnarray}
and it is not too difficult to extend it to the generic \PVI.
\vfill\eject
One such extrapolation
\begin{eqnarray}
& & {\hskip -12.0truemm}
U_r= \frac{a_{11}}{4}                                    \pmatrix{1 & 0 \cr 0 &-1 \cr}
   + \frac{(V-2X+1)\sqrt{R}}{4 X(X-1)\sqrt{V(V-1)}(V-X)} \pmatrix{0 & 1 \cr 0 & 0 \cr}
\nonumber\\ & & {\hskip -8.0truemm}
   + \frac{a_{21}}{4 X(X-1)\sqrt{R V(V-1)}}              \pmatrix{0 & 0 \cr 1 & 0 \cr},
\nonumber\\ & & {\hskip -12.0truemm}
a_{11}=-\left(\frac{1}{V}+\frac{1}{V-1}+\frac{1}{V-X}\right)\frac{\D V}{\D X}
       +\frac{\D \log R}{\D X}
			 +\frac{2}{V-X}
\nonumber\\ & & {\hskip -8.0truemm}
			 +\left(\frac{-2}{V-1}+\frac{1}{V-X}\right)F,
\nonumber\\ & & {\hskip -12.0truemm}
a_{21}=\frac{V-2X}{V-X} \left[X(X-1)\frac{\D V}{\D X}\right]^2
\nonumber\\ & & {\hskip -8.0truemm}
  +2 X(X-1)\frac{\D V}{\D X}
	\left[2 V \partial_V F -\left(\frac{2}{V-1}+\frac{1}{V-X}\right)F\right]
\nonumber\\ & & {\hskip -8.0truemm}
	+4 X(X-1)V \left(\partial_X F +\frac{F}{V-X}\right) 
  +2 X(X-1)\frac{V^2(V-1)}{V-X}
\nonumber\\ & & {\hskip -8.0truemm}
	+ V \left(-\frac{2}{V-1}+\frac{1}{V-X}\right) \left(F^2+k (V-X)^2\right)
\nonumber\\ & & {\hskip -8.0truemm}
	+2V^2(V-1)(V-X)
	 \left[\theta_\infty^2-\theta_0^2 \frac{X}{V^2}
	      +\theta_1^2 \frac{X-1}{(V-1)^2}-\theta_X^2 \frac{X(X-1)}{(V-X)^2}\right],
\nonumber\\ & & {\hskip -12.0truemm}
V_r=\frac{X(X-1)V'+F}{2(V-X)} \pmatrix{1 & 0 \cr 0 & -1 \cr}
    -\frac{(V-1)\sqrt{R}}{2 \sqrt{V(V-1)}(V-X)}  \pmatrix{0 & 1 \cr 0 & 0 \cr}
\nonumber\\ & & {\hskip -8.0truemm}
		+\frac{V\left[(X(X-1)V'+F)^2+k(V-X)^2\right]}{2 \sqrt{R V(V-1)}(V-1)} 
		\pmatrix{0 & 1 \cr 0 & 0 \cr},
\label{eqGW-X-t-codim0}
\end{eqnarray}
depends on two arbitrary functions $R$ and $F(V,X)$
and one arbitrary constant $k$,
and it reduces to (\ref{eqGW-X-t-codim2}) for the values
\begin{eqnarray}
& & {\hskip -12.0truemm}
R=R_{+} R_{-},\ F=-\frac{1+g}{2} V(V-1),\ k=-\frac{K}{4}\cdot
\end{eqnarray}
 
\vfill\eject
\section{Conclusion}

There exist at least two types of solutions 
of the Gauss-Codazzi system (\ref{eqGaussCodazziR3c})
in terms of \PVI.
In the first type, 
the set $\lbrace 1/H, 1/Q\rbrace$
is made of a harmonic function and the Hamiltonian of \PVI.
These are the solutions of Bonnet (1867) 
and Bobenko, Eitner and Kitaev (1997). 
The solution $1/Q$ harmonic defines a Bonnet surface,
while the solution $1/H$ harmonic defines a surface dual to a Bonnet surface
\cite[p.~77]{BE2000}.

\smallskip

The second type is presented for the first time in this article.
In this type,
the product $e^u H/Q$ and the logarithmic derivatives of $H$ and $Q$ 
with respect to the reduced independent variable
are homographic transforms of the same \PVI\ function,
and there is apparently no simple geometric variable which is harmonic.
The only such solution which we could find is 
(\ref{eqGC_Lie_red_g_solutionP6}),
and it requires two constraints among the four parameters of \PVI,
i.e.~the same number as the Bonnet and BEK solutions. 
The corresponding analytic surface 
is isothermic but 
is neither a Bonnet surface nor a surface dual to a Bonnet surface. 

\smallskip

Another point worth being emphasized is the power of the Lie approach.
Indeed,
among the set of infinitesimal point symmetries, 
there exists at least one algebra able to recover the two previously known \PVI\ solutions
(Bonnet 1867, 
Bobenko, Eitner and Kitaev 1997),
see
(\ref{eqgp1solq}),
(\ref{eqgm1solh}).
However, a limitation of this Lie approach is its failure
to introduce a spectral parameter in the reduced Gauss-Weingarten equations
(\ref{eqGW-X-t-codim2}).
\smallskip

In future research, we plan to remove the two constraints among the four parameters of \PVI.

\section*{Acknowledgments}

RC warmly acknowledges the generous support of both the
Centre de recherches math\'ematiques and
l'Unit\'e mixte internationale 3457 (CRM) du CNRS. 
RC was also partially supported by 
the Laboratoire de recherche conventionn\'e LRC-M\'eso.
AMG is pleased to thank Digiteo for financial support. 
We thank Hsieh Chun-Chung and Ma Hui for stimulating discussions.

\vfill\eject

\begin{appendix}

\section*{Appendix A. First integrals and reduced moving frames} 
\label{sectionReducedGW}

There are two methods to find a first integral to a third order first degree
ODE such as (\ref{eqBonnet-ODE3h}), (\ref{eqBEK-ODE3q}).

The first one relies on classical results \cite{PaiBSMF}.
If such a third order ODE admits an algebraic first integral
and if the second order ODE defined by this first integral has the Painlev\'e property,
then its degree $d$ (in the second derivative) is bounded
and can only take the values $d=1, 2, 3, 4, 6$ of Briot and Bouquet.
The first integral then has the necessary form
\begin{eqnarray}
& & K=\sum_{j=0}^d F_j(h',h,x) \left( h''\right)^j,
\end{eqnarray}
which makes it easy to compute.
In all cases handled in this paper, the degree is $d=2$.
This is probably the method used by Hazzidakis to obtain the first integral
of (\ref{eqBonnet-ODE3h}).
\medskip

The second method is to take advantage of the 
linear representation by the moving frame (Gauss-Weingarten) equations. 
Let us take the system (\ref{eqGCsystemra1a2}) as an example.
Under the reduction (\ref{eqGCReduca1a2}),
the PDE moving frame equations (\ref{eqGW-su2-R3R4}) 
\begin{eqnarray}
& &
\D \psi=U \psi \D z + V \psi \D \barz
\end{eqnarray}
is mapped to an ODE moving frame
\begin{eqnarray}
& &
\D \Psi=\tilde U_r \Psi \D \xi + \tilde V_r \Psi \D t,
\end{eqnarray}
in which the parameter $t$ 
and the transition matrix $P$ are defined as
(see \cite[\S 5.8.1]{CMBook} for the systematic derivation of $t$ and $P$),
\begin{eqnarray}
& & {\hskip -15.0 truemm}
t=\log g_1(z) + \log g_2(\barz),\
\psi=P \Psi,\ 
\label{eqNOTspectralparam}
\\ & & {\hskip -15.0 truemm}
P=\diag(({(\log g_1)'})^{1/4}({(\log g_2)'})^{-1/4},({(\log g_1)'})^{-1/4}({(\log g_2)'})^{1/4}).
\nonumber
\end{eqnarray}
These ODE moving frame equations
\begin{eqnarray}
& & {\hskip -8.0truemm}
\tilde U_r=\frac{1}{2}\pmatrix{
 \as/2 & \left((\tildeh-c e^{\displaystyle{\ad \xi/2}})/2 -\tildeq/\tildev\right) \tildev^{1/2}\cr 
         \left((\tildeh+c e^{\displaystyle{\ad \xi/2}})/2 -\tilder/\tildev\right) \tildev^{1/2} &  
-\as/2 \cr },\
\nonumber\\ & & {\hskip -8.0truemm}
\tilde V_r=\frac{1}{2}\pmatrix{
 \ad/2+\tildev'/(2 \tildev) & \left(-(\tildeh-c e^{\displaystyle{\ad \xi/2}}) -\tildeq/\tildev\right) \tildev^{1/2}\cr 
                              \left( (\tildeh+c e^{\displaystyle{\ad \xi/2}}) +\tilder/\tildev\right) \tildev^{1/2}  &
-\ad/2-\tildev'/(2 \tildev)\cr },\
\label{eqGW-xi-t}
\end{eqnarray}
then generate the first integral $K =8 \tr \tilde V_r^2$, Eq.~(\ref{eqGCfirsta1a2}). 

\vfill\eject

\section*{Appendix B. Hamiltonian and $\tau$-function of \PVI\ }
\label{AppendixPVI}

The \PVI\ function is by definition the general solution of the \PVI\ equation,
\begin{eqnarray}
& & {\hskip -20.0 truemm}
\frac{\D^2 V}{\D X^2}=
 \frac{1}{2} \left[\frac{1}{V} + \frac{1}{V-1} + \frac{1}{V-X} \right] V'^2
- \left[\frac{1}{X} + \frac{1}{X-1} + \frac{1}{V-X} \right] V'
\nonumber \\ & & {\hskip -20.0 truemm}
+ \frac{V (V-1) (V-X)}{2 X^2 (X-1)^2}
  \left[\theta_\infty^2 - \theta_0^2 \frac{X}{V^2} + \theta_1^2 \frac{X-1}{(V-1)^2}
        + (1-\theta_X^2) \frac{X (X-1)}{(V-X)^2} \right],
\label{eqPVI}
\end{eqnarray} 
which contains four equivalent singularities $V=\infty,0,1,X$.
Let us recall the main property of \PVI. 
This is the \textit{unique} function 
defined by the class of second order algebraic ordinary differential equations (ODE).
Indeed, all such ODEs which possess the Painlev\'e property
have for general solution
either an algebraic transform of \PVI\
or of one of its five degeneracies \PV, \dots, \PI,
or an elliptic function,
or the solution of some linear ODE. 

For each $\Pn(V,X)$, there exists a rational function $\Hu(V,V',X)$
whose only movable singularities are simple poles of residue unity,
thus defining a function $\tau=e^{\int \Hu\ \D X}$,
called \textit{$\tau$-function},
which has the same kind of movable singularities as 
 the general solution of a linear
ODE (namely no movable singularities at all, only movable simple zeroes).

For \PVI, there exist 24 choices for $\Hu$.
Rather than the expression built 
by Painlev\'e \cite[Eq.~(3)]{PaiCRAS1906} and Chazy \cite[$t$ page 341]{ChazyThese},
it is preferable to adopt the choice of Malmquist \cite{MalmquistP6},
\begin{eqnarray}
& & {\hskip -2.0 truemm}
\Hu_{\rm VI,M}=\frac{X(X-1) {V'}^2}{4 V(V-1)(V-X)} 
\label{eqHamVICu}
\\ & & 
+\frac{1}{4 X(X-1)}
\left[
-\theta_\infty^2 \left(V-\frac{1}{2}\right) 
- \theta_0^2            \left(\frac{X}{V}-\frac{1}{2}\right)
\right. \nonumber\\ & & \phantom{1234567890}\left.
+ \theta_1^2            \left(\frac{X-1}{V-1} -\frac{1}{2}\right)
-(\theta_X-1)^2         \left(\frac{X(X-1)}{V-X} +X -\frac{1}{2}\right)
\right].
\nonumber
\end{eqnarray}
The conjugate variables $q,p$
\begin{eqnarray}
& & {\hskip -10.0 truemm}
q=V,\
2 p= \frac{X(X-1) V'}{V(V-1)(V-X)} 
+\frac{\theta_0}{V}+\frac{\theta_1}{V-1}+\frac{\theta_X-1}{V-X},
\end{eqnarray}
then define the polynomial Hamiltonian $\Hqp(q,p,x)=\Hu(u,u',x)$ \cite{MalmquistP6} 
\begin{eqnarray}
& & \Hqp=\frac{q(q-1)(q-X)}{X(X-1)}
\left\lbrack
 p^2
 - \left(\frac{\theta_0}{q}+\frac{\theta_1}{q-1}+\frac{\theta_X-1}{q-X}\right)p
 \right.
 \nonumber\\ & & \phantom{123456789012345678901234}
 \left.
 -\frac{\theta_\infty^2-(\theta_0+\theta_1+\theta_X-1)^2}{4 q(q-1)}
 \right\rbrack.
\end{eqnarray}

The 24 $\Hu_{\rm VI}$ obey the same second order second degree ODE,
labeled (B-V) by Chazy \cite[page 340]{ChazyThese}.
After the normalization
\begin{eqnarray}
& &
Y=X(X-1)\Hu_{\rm VI}, 
\label{eqHPVInorm}
\end{eqnarray}
the ODE for $Y(X)$       
can be written either with a nice quaternary symmetry \cite {Okamoto1980-II},
\begin{eqnarray}
& &  {\hskip -20.0 truemm}
-Y'\left(X(X-1) Y''\right)^2 
-\left[{Y'}^2 -2 Y'(X Y'-Y)+n_1 n_2 n_3 n_4\right]^2
\nonumber \\ & &  {\hskip -20.0 truemm} \phantom{0123456789012345678}
+(Y' + n_1^2)(Y'+ n_2^2)(Y'+ n_3^2)(Y'+ n_4^2)=0,
\label{eqODEHPVI} 
\end{eqnarray}
or as a simplified equation with four complementary terms \cite[Eq.~(5.4)]{CosScou},
\begin{eqnarray}
& & {\hskip -8.0truemm} \hbox{SD-Ia}:\
-\left((X(X-1) Y''\right)^2 -4 Y' (X Y'-Y)^2 + 4 {Y'}^2 (X Y'-Y)
\nonumber
\\
& & {\hskip -8.0truemm} \phantom{12345678}
+A_0 {Y'}^2 + A_2 (X Y'-Y) + \left(A_3 + \frac{A_0^2}{4}\right) Y'+A_4=0,
\label{eqSDIa}
\end{eqnarray}
with the correspondence, 
\begin{eqnarray}
& &
\left\lbrace
\begin{array}{ll}
\displaystyle{
\Theta_X= \theta_X-1,
}\\ \displaystyle{
2 n_1=\theta_\infty-\Theta_X,\
2 n_2=\theta_\infty+\Theta_X,\
2 n_3=\theta_1-\theta_0,\
2 n_4=\theta_1+\theta_0,
}\\ \displaystyle{
2 A_0= \theta_\infty^2 + \theta_0^2 + \theta_1^2 + \Theta_X^2 
     =2(n_1^2+n_2^2+n_3^2+n_4^2),
}\\ \displaystyle{
4 A_2= -(\theta_\infty^2 - \Theta_X^2)(\theta_0^2 - \theta_1^2)
     = 16 n_1 n_2 n_3 n_4,
}\\ \displaystyle{
4 A_3=  (\theta_\infty^2 - \theta_1^2)(\theta_0^2 - \Theta_X^2)
}\\ \displaystyle{ \phantom{4 A_2}
     =-(n_1+n_2+n_3+n_4)(n_1+n_2-n_3-n_4)
}\\ \displaystyle{ \phantom{4 A_2xx} \times
		(n_1-n_2+n_3-n_4)(n_1-n_2-n_3+n_4),
}\\ \displaystyle{
32 A_4=
  (\theta_\infty^2 + \Theta_X^2)   (\theta_0^2 - \theta_1^2)^2
+ (\theta_\infty^2 - \Theta_X^2)^2 (\theta_0^2 + \theta_1^2)
}\\ \displaystyle{ \phantom{32 A_2}
      =32 (n_1 n_2 n_3 n_4)^2 (n_1^{-2}+n_2^{-2}+n_3^{-2}+n_4^{-2}). 
}
\end{array}
\right.
\label{SDIanjAjtj}
\end{eqnarray}
For reference, the inverse transformation is \cite[Table R]{Okamoto1980-II}, 
\begin{eqnarray}
& & {\hskip -9.0 truemm}
V=X+\frac{ \Theta_X X(X-1) Y''}{(Y'+(\theta_\infty+\Theta_X)^2/4)(Y'+(\theta_\infty-\Theta_X)^2/4)}
\nonumber\\ & & \phantom{1} 
+\frac{r_+}{Y'+(\theta_\infty+\Theta_X)^2/4}
+\frac{r_-}{Y'+(\theta_\infty-\Theta_X)^2/4},
\nonumber\\ & & {\hskip -9.0 truemm}
r_+ + r_-=
-Y-\frac{\theta_\infty^2+3\Theta_X^2}{8} (2 X-1)+\frac{\theta_1^2-\theta_0^2}{8},
\nonumber\\ & & {\hskip -9.0 truemm}
r_+ + r_-=\left( 
-Y-\frac{3\theta_\infty^2+\Theta_X^2}{8} (2 X-1)-\frac{\theta_1^2-\theta_0^2}{8}
\right) \frac{\Theta_X}{\theta_\infty}\cdot
\end{eqnarray}

\end{appendix}

\vfill\eject
%

\vfill\eject

\end{document}